%% file: note0.tex
\documentclass[12pt]{article}

\input{epsf}

\newcommand{\figWidth}{6.2in}

\input{noteDecl}

\title{Solving the sextic by iteration:\\
 A complex dynamical approach}

\author{
Scott Crass and Peter Doyle\\
Departments of Mathematics\\
Buffalo State College and University of California, San Diego\\
crass@buffalostate.edu\\
doyle@math.ucsd.edu}

\begin{document}

\maketitle

\input{note1}
\input{note2}

\input{note3}
\appendix
\input{noteApp}

\include{noteRef}

\end{document}

%% file: noteDecl.tex
\newtheorem{fact}{Fact}
\newtheorem{conj}{Conjecture}
\newtheorem{thm}{Theorem}

\newcommand{\CC}[1]{\ensuremath{\mathbf{C}^{#1}}}
\newcommand{\CP}[1]{\ensuremath{\mathbf{CP}^{#1}}}

\newcommand{\R}[1]{\ensuremath{\mathbf{R}^{#1}}}
\newcommand{\RP}[1]{\ensuremath{\mathbf{RP}^{#1}}}

\newcommand{\RR}{\ensuremath{\mathcal{R}}}

\newcommand{\Cu}[1]{\ensuremath{\mathcal{C}_{#1}}}
\newcommand{\Cb}[1]{\ensuremath{\mathcal{C}_{\overline{#1}}}}

\newcommand{\CCu}[1]{\ensuremath{{C}_{#1}}}
\newcommand{\CCb}[1]{\ensuremath{{C}_{\overline{#1}}}}

\newcommand{\Orb}[1]{\ensuremath{\mathcal{O}_{#1}}}

\newcommand{\A}[1]{\ensuremath{\mathcal{A}_{#1}}}

\newcommand{\Sym}[1]{\ensuremath{\mathcal{S}_{#1}}}

\newcommand{\D}[1]{\ensuremath{\mathcal{D}_{#1}}}

\newcommand{\G}{\ensuremath{\mathcal{G}}}

\newcommand{\Tet}{\ensuremath{\mathcal{T}}}
\newcommand{\Oct}{\ensuremath{\mathcal{O}}}

\newcommand{\I}{\ensuremath{\mathcal{I}}}
\newcommand{\Iu}[1]{\ensuremath{\mathcal{I}_{#1}}}
\newcommand{\Ib}[1]{\ensuremath{\mathcal{I}_{\overline{#1}}}}

\newcommand{\V}{\ensuremath{\mathcal{V}}}
\newcommand{\Va}{\ensuremath{\mathcal{V}_{3 \cdot 360}}}
\newcommand{\Vb}{\ensuremath{\mathcal{V}_{6 \cdot 360}}}

\newcommand{\F}[1]{\ensuremath{F_{#1}}}

\newcommand{\h}[1]{\ensuremath{h_{#1}}}

\newcommand{\id}{\ensuremath{\mathrm{id}}}

\newcommand{\Bub}{\ensuremath{\mathrm{bub}}}

\newcommand{\Ub}[1]{\ensuremath{U_{\overline{#1}}}}

\newcommand{\pQ}[3]{\ensuremath{p_{\:\!\overline{#1}#2_{#3}}}}

\newcommand{\Bar}[1]{\ensuremath{\overline{#1}}}

%% file: note1.tex

\section{Introduction}

Recently, P.~Doyle and C.~McMullen devised an iterative solution to the fifth degree polynomial
 \cite{DM}.  At the method's core is a rational mapping $f$ of \CP{1}\ with the icosahedral
(\A{5})
symmetry of a general quintic. Moreover, this \A{5}-\emph{equivariant} posseses \emph{nice}
dynamics:  for almost any initial point $a_0 \in \CP{1}$, the sequence of iterates $f^k(a_0)$
converges to one of the periodic cycles that comprise an icosahedral orbit.\footnote{For a
geometric description see \cite[p. 163]{DM}.} .  This breaking of \A{5}-symmetry provides for a
\emph{reliable} or \emph{generally-convergent} quintic-solving algorithm:  with almost any
fifth-degree equation, associate a rational mapping that has nice dynamics and whose attractor
consists of points from which one computes a root.  

An algorithm that solves the sixth-degree equation requires a dynamical system
with \A{6} symmetry.  Since there is no action of \A{6} on \CP{1}, attention turns to higher
dimensions.  Acting on \CP{2}\ is an \A{6}-isomorphic group of projective
transformations that was found by Valentiner~\cite{Val} in the late nineteenth century.  The
present work exploits this 2-dimensional \A{6}\ soccer ball in finding a ``Valentiner-symmetric"
rational mapping of \CP{2}\  whose dynamics \emph{experimentally appear} to be nice in the
above sense---transferred to the \CP{2}\ setting.  This map provides the central feature of a
conjecturally-reliable sextic-solving algorithm analogous to that employed in the quintic case.

\section{Solving equations by iteration}

For $n \leq 4$, the symmetric groups \Sym{n}\ act faithfully on \CP{1}. 
Corresponding to each action is a map whose nice dynamics provides for algorithmic convergence
to roots of given $n$th-degree equations.  For instance, Newton's method gives an iterative
solution to quadratic polynomials, but fails for higher degree equations.  The focus here
is upon the geometric and dynamical properties of  complex projective mappings rather than
numerical estimates. 

The search for elegant complex geometry and dynamics continues into degree 5.  Here, \A{5}\ is
the appropriate group, since \Sym{5}\ fails to act on the sphere.  This reduction in the galois
group requires the extraction of the square root of a polynomial's discriminant.  Such root-taking
is itself the result of a reliable iteration, namely, Newton's method.  The Doyle-McMullen algorithm
employs a  map with icosahedral symmetry and nice dynamics.

Pressing on to the sixth-degree leads to the 2-dimensional \A{6} action of the Valentiner group
\V.  Here, the problem shifts to one of finding a nice \V-symmetric mapping of \CP{2} from
whose attractor one calculates a given sextic's root.\footnote{ The solution-procedure follows
that of the quintic algorithm.  See Section~\ref{sec:solve}.}  Providing the overall framework
is the 2-dimensional \A{6}\ analogue of the icosahedron.  

Concerning the paper's style, many of this work's results are computational in nature.  As such,
they are called ``Facts".

\section{The Valentiner group}

\subsection*{Generating the group}

In order to produce an \A{6}-isomorphic group\footnote{A detailed treatment of the group's
generation as well as its combinatorial and invariant structure will be found in \cite{group}.} $\V
\subset \mbox{PGL}_3(\CC{})$, extend a ternary icosahedral group \I\ by the order-4 element
$Q$ that generates an octahedral group \Oct\ over one of the tetrahedral subgroups $\Tet \subset
\I$; symbolically, $\V = \left<\I,Q\right>$ and $\Oct = \left<\Tet,Q\right>$.  In ``orthogonal"
coordinates, $\I=\left<Z,P\right>$ where
\[ \begin{array}{ll}

 Z=\left( \begin{array}{ccc} -1&0&0\\0&1&0\\0&0&-1 \end{array}\right)&
\hspace*{20pt}
 P=\frac{1}{2}\left(\begin{array}{ccc}
    1&\tau^{-1}&-\tau\\
      \tau^{-1}&\tau&1\\
       \tau&-1&\tau^{-1}
\end{array} \right) 

\end{array} \]
 and $\tau=\frac{1+\sqrt{5}}{2}$.   To generate \V\ over \I\ take $\rho=e^{2\pi i/3}$ and
\[ \begin{array}{c}
 
Q = \left( \begin{array}{ccc} 1&0&0\\0&0&\rho^2\\0&-\rho&0 \end{array}
\right).

\end{array} \]

\subsection*{Icosahedral conics}

Sitting inside \A{6}\ are 12 \A{5} subgroups that decompose into 2
conjugate systems of 6.  The intersection of two \A{5}\ subgroups in the \emph{same} system is
an \A{4} that corresponds to the 12 tetrahedral rotations.  For subgroups in \emph{different}
systems the result is a dihedral group \D{5}.  Corresponding to each of these \A{5}\ subgroups
are ternary icosahedral subgroups of \V.  The group \I\ described above (call it \Ib{1}) preserves
the quadratic form
\[ \CCb{1} = x_1^2 + x_2^2 + x_3^2 \]
 and hence, the conic $\Cb{1}=\{\CCb{1}=0\}$ in \CP{2}.  The 5 remaining
 ``barred" conic \emph{forms} \CCb{m}\ that are stable under the conjugates
 \Ib{k}\ ($k=2,\ldots,6$) of \Ib{1}\ arise through the action of \V.  An ``unbarred" conic form
\CCu{k}\ is given by a linear combination of the \CCb{m} that is projectively invariant under \Iu{k}.

\subsection*{Anti-holomorphic symmetry}

In analogy to the 2 sets of tetrahedra found in the icosahedron, the 2 sets of icosahedral conics
exchange themselves under a ``bar-unbar" map 
\[ \Bub: \CP{2} \longrightarrow \CP{2} \]
 of order 2 that can take the form
\[ 
 [y_1,y_2,y_3] \rightarrow
 [\overline{y_1},\overline{y_2},\overline{y_3}]. 
\]
This map fixes point-wise the \RP{2}
\[ \{[y_1,y_2,y_3]\,|\,y_1,y_2,y_3\in \R{}\}. \]
 For $T \in \V,\ \Bub \circ T \circ \Bub \in \V$ so that $\Bar{\V}_{2 \cdot 360} = \left<\V, \Bub
\right>$ extends \V.  Such an ``anti-involution" exists for each of the 36 pairs of conics
$\{\Cb{a},\Cu{b}\}$.  The corresponding \RP{2}\ shares the \D{5}\ symmetry of the pair.

\subsection*{Special orbits}

The intersection of 2 conics in one system yields 4 tetrahedral points
that are face-centers of the respective icosahedra.  Overall, these
provide 60 point \V-orbits \Orb{60}\ and \Orb{\Bar{60}}.  Across
systems, 2 conics $\{\CCb{a},\CCu{b}\}$ meet in 36 antipodal pairs
$\{\pQ{a}{b}{1},\pQ{a}{b}{2}\}$ of icosahedral vertices and result in a \V-orbit \Orb{72}.

Associated with each of the 45 elements of order 2 in \A{6}\ is an
involution in \V\ that fixes a point and a line point-wise.  All
special orbits except \Orb{72}\ lie on these ``45-lines".  At their
intersections are the orbits \Orb{36}, \Orb{45}, \Orb{60}, and
\Orb{\Bar{60}}.  Each line contains 2 distinguished points that collectively give \Orb{90}.  The
generic points on the 45-lines belong to 180 point orbits.  A rich combinatorial geometry
appears in this configuration of conics and lines.

\subsection*{Invariants}

For a finite group action \G\ on a vector space over $k$, the Molien series 
$$M(\G)=\sum_{n=0}^\infty k[x]_n^\G\,t^n $$
provides one of the basic tools of classical invariant theory.  The coefficient of $t^k$ 
specifies the dimension of the space $k[x]_n^\G$ of homogeneous polynomials of degree $n$ that
are \G-invariant.  Molien's theorem supplies the generating function for the series.
\begin{thm}[Molien]  
Given a finite group action \G\ on a vector space, the Molien series is given by
$$ M(\G) = \frac{1}{|\G|} \sum_{C_{T} \subset \G} \frac{|C_{T}|}{\det\, (I\ -\ t\,T^{-1})} $$
where $C_T$ are conjugacy classes.
\end{thm}

In the Valentiner case the space is \CC{3} while the group is a 1-to-3 lift of \V\ to a
subgroup \Va\ of $\mbox{SU}_3(\CC{})$.  Moreover, the invariants of \V\ and \Va\ are in
one-to-one correspondence. There is also a 1-to-6 extension of \V\ to a so-called unitary
reflection group \Vb\ on \CC{3} whose generators are 45 involutions \cite[pp. 278, 287]{ST}. 
The elements of \Vb\ satisfy $\det\, T = \pm 1$ while \Va\ is the determinant 1 half of \Vb.  
\begin{fact}
The Molien Series for the linear Valentiner groups are given by
\begin{eqnarray*}
M(\Va) & = & \mbox{\large{$\frac{1\, +\, t^{45}}
 {(1\, -\, t^6)\, (1\, -\, t^{12})\, (1\, -\, t^{30})}$}}\\
&=& 1\, + t^6 + 2\, t^{12} + 2\, t^{18} + \ldots + t^{45} + \ldots\\
 M(\Vb) & = & \mbox{\large{$\frac{1}{(1\, -\, t^6)\, (1\, -\, t^{12})\, (1\, -\, t^{30})}$}}\\
&=& 1 + t^6 + 2\, t^{12} + 2\, t^{18} + \ldots.
\end{eqnarray*}
\end{fact}
The generating functions indicate 3 ``basic" forms that generate the ring of \Vb-invariants.  These
occur in degrees 6, 12, and 30.  In ``bub-coordinates" where one of the \Bub-involutions is given
by
conjugation of components, the
sixth--degree invariant is
\begin{eqnarray*}
F(y)&=&\alpha\,\sum_{m=1}^6 \CCb{m}(y)^3\\
     &=&\Bar{\alpha}\,\sum_{m=1}^6 \CCu{m}(y)^3\\
     &=&10\,y_1^3 y_2^3+9\,y_1^5 y_3+9\,y_2^3 y_3-45\,y_1^2 y_2^2 y_3^2-135\,y_1 y_2 
             y_3^4+27\,y_3^6.
\end{eqnarray*}
Here, $\alpha$ is a simplifying factor.  The forms of degrees 12 and 30 arise from Hessian and
``bordered Hessian" determinants that involve invariants of  lower degree:
 \begin{eqnarray*}
\Phi(y) &=& \alpha_{\Phi}\,|H(F(y))|\\
&=& 6\,y_{1}^{11} y_{2} + \ldots + 6\,y_{1} y_{2}^{11} + \ldots + 729\,y_{3}^{12}\\
\Psi(y) &=&  \alpha_{\Psi}\,\left| \begin{array}{c|c}
H(\Phi(y))& \begin{array}{c} F_{y_{1}}\\F_{y_{2}}\\F_{y_{3}} \end{array}\\
\hline
F_{y_{1}}\ F_{y_{2}}\ F_{y_{3}}&0
\end{array} \right|\\
&=& 3\, y_1^{30} + \ldots\ + 3\, y_2^{30} + \ldots + 57395628\,y_3^{30}.
\end{eqnarray*}
The constants $\alpha_{\Phi}=-1/20250$ and $\alpha_{\Psi}=1/24300$ remove the highest
common factor among the respective terms.  Finally, the product of the 45 linear forms that
correspond to the generating involutions is a relative \Vb-invariant but an absolute \Va-invariant
and hence, a projective \V-invariant.  This degree 45 form is given by the Jacobian determinant
$$ \begin{array}{lll}

X(y)&=&\alpha_X\,\left| \begin{array}{ccc}

F_{y_{1}}&F_{y_{2}}&F_{y_{3}}\\

\Phi_{y_{1}}&\Phi_{y_{2}}&\Phi_{y_{3}}\\

\Psi_{y_{1}}&\Psi_{y_{2}}&\Psi_{y_{3}}

\end{array} \right| \\[20pt]

&=& y_1^{45} + \ldots - y_2^{45} + \ldots + 3570467226624\,y_2^5 y_3^{40}

\end{array} $$ 
where $\alpha_X=-1/4860$. Being \Vb-invariant, $X^{2}$ is a polynomial
in $F$, $\Phi$, $\Psi$: 
\begin{eqnarray} \label{eq:X^2}
X^{2} 
&=&\frac{1}{19683}\,(4\,F^{13} \Phi + 80\,F^{11} \Phi^2 + 816\,F^9 \Phi^3 + 
4376\,F^7 \Phi^4 +\\ \nonumber 
&&13084\,F^5 \Phi^5 + 12312\,F^3 \Phi^6 + 5616\,F \Phi^7 + 
18\,F^{10} \Psi +\\ \nonumber 
&&198\,F^8 \Phi
\Psi + 954\,F^6 \Phi^2 \Psi - 198\,F^4 \Phi^3 \Psi - 5508\,F^2 \Phi^4 \Psi 
-\\ \nonumber 
&&1944\,\Phi^5 \Psi - 162\,F^5 \Psi^2 -1944\,F^3 \Phi \Psi^2
- 1458\,F \Phi^2 \Psi^2 + 729 \Psi^3). 
\end{eqnarray}

%% file: note2.tex

\section{Rational maps with Valentiner symmetry}

An iterative solution of the sextic involves a family of $A_{6}$-symmetric dynamical systems
that is parametrized by sixth-degree polynomials. This parametrization amounts to conjugation of
a fixed $A_{6}$-equivariant $f$ by a map
$$S_{p}:\CP{2} \rightarrow \CP{2}$$
that is associated with a given sextic $p$ to be solved.  Accordingly,
the map $f$ is the centerpiece of a sextic-solving algorithm.\footnote{An expanded discussion of
\V-symmetric maps will appear in \cite{sextic}.}

\subsection{Finding equivariant maps}

Since \V-equivariants associate one-to-one with \Va-equivariants, the search for symmetric maps
can take place within the regime of the linear action.  The linear Valentiner groups \Va\ and \Vb\
act on the exterior algebra $\Lambda(\CC{3})$ by
$$ (T(\alpha))(x) = \alpha(T^{-1}x)  $$
where $T\in\,\Va$ or $T\in\Vb$, $\alpha$ is a 0,1,2,3-form, and $x\in\CC{3}$.

This action provides guidance in the search for equivariants.  Such
utility is due to a correspondence between \G-invariant forms and \G-equivariant maps.

\begin{thm}
 For a given finite action \G\ on \CC{3} and a $\G$-invariant 2-form
\begin{eqnarray*}
 \phi(x)&=&f_{1}(x)\, dx_2 \wedge dx_3 + f_{2}(x)\, dx_3 \wedge
dx_1 + f_{3}(x)\, dx_1 \wedge dx_2,
  \end{eqnarray*} 
the map \mbox{$f=(f_{1},f_{2},f_{3})$} is \G-equivariant.  
\end{thm}

There is a 2-variable ``exterior" Molien series $M(\Lambda^{\G})$ in which the
variables $s$ and $t$ index respectively the rank of the form and the polynomial degree:
$$ M(\Lambda^{\G}) = \sum_{p=0}^{n} \left(\,
  \sum_{m=0}^{\infty} \left(\dim\,\Lambda_p \left(\CC{n} \right)
  _m^\G \right)
    t^{m} \right) s^{p} $$
where $\Lambda_p(\CC{n})_m^\G$ are the $\G$-invariant $p$-forms of degree $m$ (\cite[p.
62]{Benson} or \cite[pp. 265ff]{Smith}).  Projection of $\Lambda_p(\CC{n})_m$ onto
$\Lambda_p(\CC{n})^\G_m$ yields the analogue to Molien's theorem.
\begin{thm}
The exterior Molien series for a finite group action \G\ is given by the generating function:
 \begin{eqnarray*}
M(\Lambda^{\G})&=&
\frac{1}{\left|\G\right|}\
    \sum_{\mathcal{C}_{T} \subset \G}
     \left|\mathcal{C}_{T} \right|\,
      \frac{\det\left(I + s\,T^{-1}\right)}
           {\det\left(I - t\,T^{-1}\right)}
\end{eqnarray*}
 where $\mathcal{C}_{T}$ are conjugacy classes.  
\end{thm}
An expansion of the 2-variable series $M(\Lambda^{\Va})$ will appear  in \cite{sextic}.

In the series $M(\Lambda^{\mathcal{H}})$ for a subgroup $\mathcal{H}
\subset\:\G$ the contribution of the \G-invariant 0-forms $\Lambda_0^\G$ disappears through
division of $M(\Lambda^{\mathcal{H}})$ by $M(\Lambda_0^\G)=M(\CC{}[x]^\G)$.
The resulting \emph{polynomial} in 2 variables displays the degrees of the generating
$\mathcal{H}$-invariant forms.  

\begin{fact}
For the Valentiner group the ``exterior Molien quotient"  is
$$\begin{array}{ccl}

M(\Lambda^{\Va})/ M(\Lambda_{0}^{\Vb})&=&
(1+t^{45})\,+\,(t^{5}+t^{11}+t^{20}+t^{26}+t^{29}+t^{44})s\,+\\

&&(t+t^{16}+t^{19}+t^{25}+t^{34}+t^{40})s^{2}\,+\,(1+t^{45})s^{3}.  

\end{array}$$
\end{fact}

Given \G-invariants $G(x)$ and $H(x)$, the 2-form
$$dG(x) \wedge dH(x) = (\nabla G(x) \times \nabla H(x)) \cdot (dx_2 \wedge dx_3,dx_3 \wedge  
  dx_1,dx_1 \wedge dx_2) $$
is also \G-invariant.   Here, `$\nabla$' and `$\cdot$' stand for a ``formal" gradient $\nabla F =
(\frac{\partial
F}{\partial y_1},\frac{\partial F}{\partial y_2},\frac{\partial
F}{\partial y_3})$ and dot product while `$\times$' specifies the cross-product.  Thus, exterior
algebra accounts for
\V-equivariants of degrees 16, 34, and 40:
\begin{eqnarray*}
  \psi_{16}(y)&=&\nabla F(y) \times \nabla \Phi(y)\\
  \phi_{34}(y)&=&\nabla F(y) \times \nabla \Psi(y)\\
     f_{40}(y)&=&\nabla \Phi(y) \times \nabla \Psi(y).
\end{eqnarray*}

\subsection*{A family of degree 19 maps}

The Molien series for Valentiner equivariants
$$ t+t^{7}+2\,t^{13}+t^{16}+3\,t^{19}\ldots $$
specifies 3 dimensions worth of maps in degree 19 of which 2 are due
to promotion of the identity by the 2 dimensions of degree 18
invariants.  Hence, there are, as the exterior Molien quotient indicates, non-trivial
\V-symmetric maps in degree 19.  How do these arise?  Since there is no
apparent exterior algebraic means of producing such a map, the more
practical matter of computing them takes priority.

\subsubsection*{19 = 64 - 45}

Multiplication of a degree 19 equivariant
$f$ by $X_{45}$ yields $X\!\cdot\! f$ in the 14 dimensional space of
64-maps.  There are 14 ways of promoting the maps $\psi_{16},\phi_{34},f_{40}$ to degree 64:

\begin{description}

\item{1)} 7 dimensions of degree 48 invariants to promote $\psi_{16}=\nabla F \times \nabla
\Phi$

\item{2)} 4 dimensions of degree 30 invariants to promote $\phi_{34}=\nabla F \times \nabla
\Psi$
  
\item{3)} 3 dimensions of degree 24 invariants to promote $f_{40}=\nabla \Phi \times \nabla
\Psi$.

\end{description}
These 14 maps span the space of degree 64 equivariants.  Thus, $X\!\cdot\! f$
is a combination of maps whose computation is straightforward.
Reasoning in the other direction, a 64-map
$$ 
f_{64} = \alpha\,F_{48}\cdot\psi_{16} + \beta\,F_{30}\cdot\phi_{34} +   
 \gamma\,F_{24}\cdot f_{40}
$$
 that ``vanishes" on the 45-lines---i.e., each coordinate function of
$f_{64}$ vanishes---must have a factor of $X$ which can then be divided
away.  The quotient is a degree 19 equivariant
$$f_{19}=\frac{f_{64}}{X_{45}}.$$

To arrange for the vanishing of $f_{64}$ on the 45-lines, only one line need
be considered;  symmetry tends to the remaining lines.  Forcing
$f_{64}$ to vanish at 12 ``independent" points on a
45-line\footnote{The points are \emph{independent} in the sense that
the 12 resulting linear conditions in the 14 undetermined coefficients
of $f_{64}$ are independent.} yields a 2-parameter family of 64-maps
each member of which vanishes on $\{X=0\}.$  The 2 \emph{inhomogeneous}
parameters reflect the 3 dimensions (i.e., homogeneous parameters)
worth of degree 19 \V-equivariants. In \Bub-coordinates, setting
these 2 parameters equal to 0 and normalizing the coefficients yields the map:
\begin{eqnarray} \label{eq:f64}
 f_{64}(y) &=&\left[10\,F(y)^6\,\Phi(y) + 100\,F(y)^4\,\Phi(y)^2 +
  45\,F(y)^2\,\Phi(y)^3 + 156\,\Phi(y)^4\ \right. + \nonumber\\
&&\left. 39\,F(y)^3\,\Psi(y) + 51\,F(y) \Phi(y)\,\Psi(y)\right]\cdot
   \psi(y) - \nonumber\\
&& 27\,\Psi(y) \cdot \phi(y) + 54\,\Phi(y)^2 \cdot f(y).
\end{eqnarray}
 The 2-parameter family of non-trivial 19-maps is then given by
\begin{eqnarray}  \label{eq:19-maps}
 g_{19}(y;a,b)\ =\ f_{19}(y)\:+\:  \left(a\,F(y)^{3} +
  b\,F(y)\,\Phi(y)\right) \cdot y.
\end{eqnarray}
 Are any of these maps dynamically ``special"?  Indeed, what might
it mean to be special in this sense?

\subsubsection*{Extended symmetry in degree 19}

Since $f_{19}$ is a non-trivial $\Bar{\V}_{2 \cdot
360}$-equivariant---note the integer coefficients in (\ref{eq:f64}), each member of the following
1-parameter family, being impartial towards the 2 systems of conics, also enjoys the additional
symmetry:
$$ f_{19}\ +\ a\,F\left(B_{12}+U_{12}\right)\cdot \id $$
 where $B_{12} = \prod C_{\overline{k}}$ and $U_{12} = \prod C_{k}$ are
the degree 12 invariants given by the product of the respective 6 conic
forms.  To honor the doubled symmetry a member of this family must set-wise fix
each \RP{2}\ among the 36 that are point-wise fixed by the order 2 \Bub\ maps.

\section{The degree 19 map}   \label{sec:19map}

Since 19 is an ``equivariant degree" for the  binary icosahedral group \cite[p. 166]{DM}, there
arises the prospect of finding a \V-equivariant that restricts to self-mappings of the
conics.\footnote{The forthcoming \cite{sextic} will give details of  the procedure used in finding
this special map.}  Appropriate expenditure of the 2 parameters in degree 19 purchases such a
map.
\begin{fact}
There is a unique degree 19 \V-equivariant $h_{19}$ that preserves all 12 of the icosahedral
conics.
\end{fact}
 On each conic-icosahedron the map has the following geometric description:  stretch each face
\textsf{F} over the 19 faces in the complement of the face antipodal to \textsf{F} while making a
half-turn in order to send the 3 vertices and edges of \textsf{F} to their antipodal vertices and
edges \cite[p. 163]{DM}.  The 20 face-centers are fixed and repelling, while the 12 period-2
vertices form the map's critical set.  Consequently \cite[p. 156]{DM}, the map has nice dynamics.
\begin{thm}
Under $h_{19}$, the trajectory of almost any point on an icosahedral conic tends to an antipodal
pair of the superattracting vertices.
\end{thm}
In ``bub-coordinates", the conic-fixing equivariant has the expression:
 \vspace{10pt}

\input{h19}

\subsection*{Dynamical behavior}

The discovery of $h_{19}$ supplies the unique degree 19 equivariant that
self-maps, in addition to the 45-lines and the 36 \Bub-\RP{2}s, the 12
conics.  The dynamics \emph{on} each conic is well-understood:  the critical icosahedral vertices
attract a full measure's worth of points.\footnote{See Figure~\ref{fig:h19Conic} for a basins
of attraction plot.}  Indeed, the map attracts on a full \CP{2}\
neighborhood of such a pair of 72-points. Hence, the Fatou components of the restricted maps
$h_{19}|_{\Cb{a}}$ are the intersections with \Cb{a} of the map's Fatou components in \CP{2}. 
Is this attracting behavior of the conics pervasive in the measure-theoretic sense?   What about the
``restricted" dynamics on the 45-lines and $\RP{2}$s?  On a \Bub-\RP{2} the experimental
evidence\footnote{See Figure~\ref{fig:h19RP2} for a plot of the attracting basins on one such
\RP{2}.} strongly suggests that 
\begin{itemize}
\item[1)]  the 72-points are the only attractors, 
\item[2)]  the set of 45-lines $\{X=0\}$ is repelling, and 
\item[3)]  there is no region with thickness or, indeed, positive measure that remains outside of
their influence.  
\end{itemize}
What significance does the \RP{2}-dynamics hold for the \CP{2}-dynamics?  Extensive
trials\footnote{At \emph{http://math.ucsd.edu/\~{}scrass} there
are \emph{Mathematica} notebooks and supporting files with which to iterate $h_{19}$ and,
from the output, to approximate a solution to a sixth-degree equation.} on \CP{2}\ have not
revealed behavior contrary to that observed on the \RP{2}s.
\begin{conj}  \label{conj:h19}
The only attracting periodic points for $h_{19}$ are the elements of \Orb{72}.  Moreover, the
union of the basins of attraction for \Orb{72}\ has full measure in  \CP{2}.
\end{conj}

The 45-lines present a problem in that they map to themselves but do not
contain the 72-points.\footnote{Again, this is a feature peculiar to
the 72-points.  They form the only special \V-orbit that does not lie
on the 45-lines.} The basin plot in the Appendix reveals repelling behavior along the
\RP{1}\ where the \RP{2}\ meets one of the 45-lines that are set-wise fixed by
the \RP{2}'s \Bub map.  
\begin{conj}
On the 45-lines \h{19}\ is repelling and, hence, $\{X=0\}$ resides in the Julia set $J_{h_{19}}$.
\end{conj}
Is $J_{h_{19}}$ the closure of the backward orbit of the 45-lines?

%% file: h19.tex

\noindent $h_{19}(y) = 1620\,F(y)^3 \cdot [y_1,y_2,y_3] + f_{19}(y) = $
\vspace{10pt}
\begin{sloppypar} \raggedright \noindent
$[-3591\,y_1^{15} y_2^4 - 
   5263\,y_1^{10} y_2^9 + 
   9747\,y_1^5 y_2^{14} - 81\,y_2^{19} + 
   17955\,y_1^{12} y_2^6 y_3 +  
   10260\,y_1^7 y_2^{11} y_3 - 
   7695\,y_1^2 y_2^{16} y_3 - 
   107730\,y_1^{14} y_2^3 y_3^2 - 
   74385\,y_1^9 y_2^8 y_3^2 +  
   161595\,y_1^4 y_2^{13} y_3^2 - 
   969570\,y_1^{11} y_2^5 y_3^3 + 
   1292760\,y_1^6 y_2^{10} y_3^3 - 
   46170\,y_1 y_2^{15} y_3^3 -   
   2346975\,y_1^8 y_2^7 y_3^4 - 
   807975\,y_1^3 y_2^{12} y_3^4 - 
   3587409\,y_1^{10} y_2^4 y_3^5 + 
   10277442\,y_1^5 y_2^9 y_3^5 +   
   13851\,y_2^{14} y_3^5 - 
   969570\,y_1^{12} y_2 y_3^6 - 
   3986010\,y_1^7 y_2^6 y_3^6 - 
   1939140\,y_1^2 y_2^{11} y_3^6 -  
   5263380\,y_1^9 y_2^3 y_3^7 - 
   28117530\,y_1^4 y_2^8 y_3^7 + 
   831060\,y_1^{11} y_3^8 + 
   2423925\,y_1^6 y_2^5 y_3^8 +   
   4363065\,y_1 y_2^{10} y_3^8 - 
   24931800\,y_1^8 y_2^2 y_3^9 + 
   43630650\,y_1^3 y_2^7 y_3^9 -
   31123197\,y_1^5 y_2^4 y_3^{10}+   
   9598743\,y_2^9 y_3^{10} + 
   14959080\,y_1^7 y_2 y_3^{11} +
   23269680\,y_1^4 y_2^3 y_3^{12} - 
   26178390\,y_1^6 y_3^{13}+   
   52356780\,y_1 y_2^5 y_3^{13} + 
   18698850\,y_1^3 y_2^2 y_3^{14} + 
   20194758\,y_2^4 y_3^{15} + 
   22438620\,y_1^2 y_2 y_3^{16}+   
   7479540\,y_1 y_3^{18},$\\[5pt] 
$-81\,y_1^{19} + 9747\,y_1^{14} y_2^5 - 
   5263\,y_1^9 y_2^{10} - 
   3591\,y_1^4 y_2^{15} -   
   7695\,y_1^{16} y_2^2 y_3 + 
   10260\,y_1^{11} y_2^7 y_3 + 
   17955\,y_1^6 y_2^{12} y_3 + 
   161595\,y_1^{13} y_2^4 y_3^2 -   
   74385\,y_1^8 y_2^9 y_3^2 - 
   107730\,y_1^3 y_2^{14} y_3^2 - 
   46170\,y_1^{15} y_2 y_3^3 +
   1292760\,y_1^{10} y_2^6 y_3^3-   
   969570\,y_1^5 y_2^{11} y_3^3 - 
   807975\,y_1^{12} y_2^3 y_3^4 -
   2346975\,y_1^7 y_2^8 y_3^4 + 
   13851\,y_1^{14} y_3^5+   
   10277442\,y_1^9 y_2^5 y_3^5 - 
   3587409\,y_1^4 y_2^{10} y_3^5 - 
   1939140\,y_1^{11} y_2^2 y_3^6 - 
   3986010\,y_1^6 y_2^7 y_3^6 -   
   969570\,y_1 y_2^{12} y_3^6 - 
   28117530\,y_1^8 y_2^4 y_3^7 - 
   5263380\,y_1^3 y_2^9 y_3^7 + 
   4363065\,y_1^{10} y_2 y_3^8+   
   2423925\,y_1^5 y_2^6 y_3^8 + 
   831060\,y_2^{11} y_3^8 + 
   43630650\,y_1^7 y_2^3 y_3^9 - 
   24931800\,y_1^2 y_2^8 y_3^9 +
   9598743\,y_1^9 y_3^{10} -
   31123197\,y_1^4 y_2^5 y_3^{10} + 
   14959080\,y_1 y_2^7 y_3^{11} + 
   23269680\,y_1^3 y_2^4 y_3^{12} +   
   52356780\,y_1^5 y_2 y_3^{13} - 
   26178390\,y_2^6 y_3^{13} + 
   18698850\,y_1^2 y_2^3 y_3^{14} + 
   20194758\,y_1^4 y_3^{15}+   
   22438620\,y_1 y_2^2 y_3^{16} + 
   7479540\,y_2 y_3^{18},$\\[5pt]
$-1026\,y_1^{17} y_2^2 - 
   3078\,y_1^{12} y_2^7 - 
   3078\,y_1^7 y_2^{12} - 
   1026\,y_1^2 y_2^{17}-   
   5130\,y_1^{14} y_2^4 y_3 + 
   113240\,y_1^9 y_2^9 y_3 - 
   5130\,y_1^4 y_2^{14} y_3 + 
   3078\,y_1^{16} y_2 y_3^2 -  
   272916\,y_1^{11} y_2^6 y_3^2 - 
   272916\,y_1^6 y_2^{11} y_3^2 + 
   3078\,y_1 y_2^{16} y_3^2 + 
   215460\,y_1^{13} y_2^3 y_3^3 +   
   687420\,y_1^8 y_2^8 y_3^3 + 
   215460\,y_1^3 y_2^{13} y_3^3 + 
   4617\,y_1^{15} y_3^4 + 
   937251\,y_1^{10} y_2^5 y_3^4 +   
   937251\,y_1^5 y_2^{10} y_3^4 + 
   4617\,y_2^{15} y_3^4 + 
   290871\,y_1^{12} y_2^2 y_3^5 + 
   4813992\,y_1^7 y_2^7 y_3^5 +   
   290871\,y_1^2 y_2^{12} y_3^5 - 
   1454355\,y_1^9 y_2^4 y_3^6 - 
   1454355\,y_1^4 y_2^9 y_3^6 + 
   2520882\,y_1^{11} y_2 y_3^7 +   
   8812314\,y_1^6 y_2^6 y_3^7 + 
   2520882\,y_1 y_2^{11} y_3^7 + 
   19876185\,y_1^8 y_2^3 y_3^8 + 
   19876185\,y_1^3 y_2^8 y_3^8 -  
   2036097\,y_1^{10} y_3^9 + 
   5623506\,y_1^5 y_2^5 y_3^9 - 
   2036097\,y_2^{10} y_3^9 + 
   5235678\,y_1^7 y_2^2 y_3^{10} +   
   5235678\,y_1^2 y_2^7 y_3^{10} + 
   37813230\,y_1^4 y_2^4 y_3^{11} - 
   2617839\,y_1^6 y_2 y_3^{12} - 
   2617839\,y_1 y_2^6 y_3^{12} -   
   2908710\,y_1^3 y_2^3 y_3^{13} + 
   6357609\,y_1^5 y_3^{14} + 
   6357609\,y_2^5 y_3^{14} - 
   5983632\,y_1^2 y_2^2 y_3^{15} - 
   4487724\,y_1 y_2 y_3^{17} - 
   1023516\,y_3^{19}]$.
\end{sloppypar}

%% file: note3.tex

\section{Solving the sextic}    \label{sec:solve}

Consider the rational functions:
$$Y_1 = \alpha\, \frac{\Phi}{F^2}\hspace*{25pt}Y_2 = \beta\,\frac{\Psi}{F^5} $$
where $\alpha$ and $\beta$ are chosen so that $Y_1$ and $Y_2$ each
 take on the value 1 at a 36-point.\footnote{In the bub-coordinates
 $y$, $\alpha = 1$ and $\beta = 1/4$.}  Given values $a_1$ and $a_2$ of
 $Y_1$ and $Y_2$, the ``form-problem" is to find a point $z$ in \CP{2}\ that
 belongs to the \V-orbit $Y_1^{-1}(a_1) \cap Y_2^{-1}(a_2)$. Solving the Valentiner
form-problem is tantamount to solving  the sextic \cite[pp. 308-10]{Fr}.  

By algebraic transformation, the general 6-parameter sextic $p(x)$ can be reduced to a
``resolvent" that depends on the 2 parameters $Y_1$ and $Y_2$.  In addition to the square root
of $p$'s discriminant, such a reduction requires the extraction of a cube root---a so-called
``accessory irrationality", since its adjunction to the coefficient field does not reduce the galois
group of $p$ \cite[p. 285]{Fr}.

\begin{fact}
The parameters $Y_1$ and $Y_2$ specify a family of sextic Valentiner resolvents
$$R_z(u) = \prod_{m=1}^{6}\left(u - U_{\Bar{m}}\right) $$
where the roots are
$$U_{\Bar{m}} = \frac{\Cb{m}^3(z)}{F(z)}\,,\ m=1,\ldots,6. $$
\end{fact}
These resolvents are especially suited for solution by an iterative
algorithm\footnote{The algorithm sketched below will receive full treatment in
\cite{sextic}.} that exploits Valentiner symmetry and symmetry-breaking.  
\begin{fact}
Each resolvent $R_z(u)$ is \V-invariant in the parameter z.  Accordingly, $R_z(u)$ is expressible
in $Y=(Y_1,Y_2)$:

$$ \begin{array}{lll}

R_Y(u) &=& u^{6}\:+\:\frac{-5\,+\,\sqrt{15}\,i}{90}\, u^5\,+\\[5pt]

&&\frac{11\,(1\,-\,\sqrt{15}\,i)\,-\,
 3\,(3\,+\,\sqrt{15}\,i)\,Y_1}{2^2 3^5 5^2}\, u^4
  \:+\:\frac{(100\,+\,57\,\sqrt{15}\,i)\,+\,
   9\,(30\,+\,\sqrt{15}\,i)\,Y_1}{3^9 5^4}\, u^3\,+\\[5pt]

&&\frac{-(152\,+\,17\,\sqrt{15}\,i)\,+\,
 18\,(-21\,+\,4\,\sqrt{15}\,i)\,Y_1\,+\,
  27\,(-4\,+\,\sqrt{15}\,i)\,Y_1^2}{2^2 3^{11} 5^5}\,u^2\,+\\[5pt]

&&\frac{(425\,+\,103\,\sqrt{15}\,i)\,+\,
 6\,(75\,+\,193\,\sqrt{15}\,i)\,Y_1\,+\,
  27\,(-25\,+\,33\sqrt{15}\,i)\,Y_1^2\,-\,7776\,\sqrt{15}\,i\,Y_2}
   {2^3 3^{14} 5^8}\,u\,+\\[5pt]

&&\frac{-(5\,+\,3\,\sqrt{15}\,i)\,+\,
 9\,(15\,-\,7\,\sqrt{15}\,i)\,Y_1\,+\,
  81\,(25\,-\,\sqrt{15}\,i)\,Y_1^2\,+\,
   81\,(45\,+\,11\sqrt{15}\,i)\,Y_1^3}{2^4 3^{18} 5^8}\,.

\end{array} $$
\end{fact}

\subsection*{Parametrized families of Valentiner groups}

The algorithm that solves a given resolvent $R_Y$ employs an iteration
of a dynamical system that belongs to a family of degree-19 maps $\h{Y}(w)$
that are 1) conjugate versions of $h_{19}(y)$ and 2) parametrized by $Y={(Y_1,Y_2)}$.

Consider the family of projective transformations in $\{w_1,w_2,w_3\}$ that are degree 25 in the
parameters $\{z_1,z_2,z_3\}$: 
$$ y = \tau_z(w) = [F(z)^4\cdot z]\,w_1\,+\,[F(z)\cdot h_{19}(z)]\,
  w_2\, +\, k_{25}(z)\, w_3. $$
Here, $k_{25}(z)$ is the equivariant that in ``Hermitian\footnote{The elements $T
\in \Va$ satisfy $T \Bar{T^{t}} = I$.} coordinates" $v$ is given by
$$ k_{25}(v) = \Bar{\nabla F(\Bar{\nabla F(v)})}. $$
Much of this development amounts to keeping track of coordinates.  The
Valentiner actions $\V_z$ and $\V_y$ on the respective planes
$\CP{2}_z$ and $\CP{2}_y$ are the same with $z$ merely replacing $y$.
Meanwhile, the group $\V_w=\tau_z^{-1}V_y\tau_z$ acts on the $w$-plane
$\CP{2}_w$.  

By invariance in $z$, $F(\tau_z(w))$ admits expression in $Y$:
$$ F_Y(w) = \frac{F(\tau_z(w))}{\alpha\, F(z)^{25}} $$
where $\alpha$ is the highest common factor over the coefficients.  As for the remaining basic
invariants, they arise from the sixth-degree form as before.  However, a parametrized
change of coordinates requires special handling:
\begin{eqnarray*} 
\Phi(y)&=&\frac{(\alpha\,F(z)^{25})^3}{|\tau_z|^2}\,\Phi_Y(w)\\
\Psi(y)&=&\frac{(\alpha\,F(z)^{25})^8}{|\tau_z|^6}\,\Psi_Y(w)\\
 X(y)&=&\frac{(\alpha\,F(z)^{25})^{12}}{|\tau_z|^9}\,X_Y(w).
\end{eqnarray*}
With an invariant system in place for each (non-singular) value of $Y$, computation of
the degree 19 map $h_Y(w)$ that preserves all 12 of the conics follows the original procedure.

\subsection*{Root-finding}

Under Conjecture~\ref{conj:h19}, the trajectory $\{h_{19}^k(y_0)\}$ converges
to a pair of 72-points for almost any $y_0 \in \CP{2}_y$.
Being conjugate to $h_{19}(y)$, the maps $h_Y(w)$ share this property
for points in $\CP{2}_w$. By breaking the \A{6}\ symmetry of $R_Y(u)=0$,
$h_Y(w)$ qualifies for a role in root-finding.

Now, consider the rational function
$$ \Bar{J}_z(w) = \frac{\Bar{\Gamma}_z(w)^3}{F(z)\,\Psi(\tau_z(w))} $$
where
$$ 
\Bar{\Gamma}_z(w) = \sum_{m=1}^6 \prod_{n \neq m}
\CCb{n}(\tau_z(w))\cdot\CCb{m}(z). 
$$ 
 At a 72-point pair $\{q_1,q_2\}=\tau_z^{-1}(\{\pQ{a}{b}{1},\pQ{a}{b}{2}\})$ in $w$-space
five of the six terms in $\Bar{\Gamma}_z(w)$ vanish.  The result is a ``selection" of one of the
six roots \Ub{a}(z) of $R_Y(u)$:
\begin{eqnarray*} 
 \Bar{J}_z(\{q_1,q_2\})
&=& 
  \frac{\left[\prod_{n \neq a} \CCb{n}(\{\pQ{a}{b}{1},\pQ{a}{b}{2}\})\right]^3}
       {\Psi(\{\pQ{a}{b}{1},\pQ{a}{b}{2}\})}\,\frac{\CCb{a}(z)^3}{F(z)}\\
&=&\frac{\Ub{a}(z)}{\mu}\,.
\end{eqnarray*}
 Here, $\mu$ is the value of  
$$ 
\frac{\Psi(\{\pQ{a}{b}{1},\pQ{a}{b}{2}\})}
        {[\prod_{n \neq a} \CCb{n}(\{\pQ{a}{b}{1},\pQ{a}{b}{2}\})]^3}\
$$
which is constant on \Orb{72}:
$$ \mu = \frac{6561\,(279 + 145 \sqrt{15}\,i)}{2}. $$

In light of their $\V_z$-invariance, $\Bar{\Gamma}_z(w)$ and $\Bar{J}_z(w)$
are expressible in $Y$ and $w$:
$$ \Bar{\Gamma}_Y(w) = \frac{\Bar{\Gamma}_z(w)}{\beta\,F(z)^{42}} $$
where $\beta$ is a simplifying factor.  As for $\Bar{J}_z(w)$,
let $T_Y$ be the polynomial in $Y_1$ and $Y_2$ that satisfies
$$ |\tau_z|^2 = F(z)^{25}\,T_Y $$
so that
\begin{eqnarray*}
\Psi(\tau_z(w))&=&\frac{\alpha^8\, F(z)^{125}}{T_Y^3}\,\Psi_Y(w)
\end{eqnarray*}
 and finally,
$$ \Bar{J}_Y(w) = \frac{\beta^3}{\alpha^8}\,
            \frac{(T_Y\,\Bar{\Gamma}_Y(w))^3}{\Psi_Y(w)}\,.
$$

To summarize the procedure:

\begin{itemize}

\item[1)]  Select a value $A=(A_1,A_2)$ of $Y=(Y_1,Y_2)$ and, thereby,
a sixth-degree resolvent $R_A(u)$. (For sake of description, let $z \in
Y_1^{-1}(A_1)\cap Y_2^{-1}(A_2)$.  The algorithm actually finds a root without explicitly
inverting $Y_1$ or $Y_2$.)

\item[2)]  From an initial point $p \in \CP{2}_w$, iterate the map $h_A(w)$ to convergence:
$$ h^n_A(p) \longrightarrow \{q_1,q_2\} \in (\Orb{72})_w \subset \CP{2}_w. $$
As output take the pair of approximate 72-points in $\CP{2}_w$
$$\{p_1,p_2\}\approx \{q_1,q_2\} =
\{\tau^{-1}_z(\pQ{a}{b}{1}),\tau^{-1}_z(\pQ{a}{b}{2})\}. $$

\item[3)]  Approximate a root of $R_A(u)$:
\begin{eqnarray*}
\Ub{a}(z)&\approx&\mu\,\Bar{J}_A(p_1)\\
               &\approx&\frac{\mu\,\beta^3\,(T_A\,\Bar{\Gamma}_A(p_1))^3}
                  {\alpha^8\,\Psi_A(p_1)}\,.
\end{eqnarray*}

\end{itemize}

%% file: noteApp.tex

\section{Evidence for Conjecture \ref{conj:h19}} 

\label{sec-Seeing}

The following pictures provide empirical dynamical information for the
special map $h_{19}$. The program \emph{Dynamics}, running on a Silicon
Graphics Indigo-2, created the basin plots with the ``BAS"
routine.\footnote{The software is the work of H. Nusse and J. Yorke
while E. Kostelich is responsible for the Unix implementation. See their
manual \cite{NY}.} This procedure colors a grid-cell in the event that
the trajectory of the cell's center gets close enough to a 72-point to
guarantee ultimate attraction to the associated period-2 cycle. The
color depends upon the destination of its center. All plots have the
maximum resolution available: a \mbox{$720 \times 720$} grid of cells.
If, in a specified number of iterations, the center's trajectory fails
to converge to a pair of 72-points, the cell is left uncolored. Such
points are candidates for the Julia set. According to
Conjecture~\ref{conj:h19}, the basins of the 72-points fill up \CP{2}\
in the measure theoretic sense. The pictures do not belie the claim.

In the basin plot for $h_{19}$ restricted to an affine plane in one of
the \Bub-\RP{2}s \RR, the chosen coordinates make evident the map's
$D_{5}$ symmetry; the 1-point orbit lies at the origin while its
companion 1-line orbit resides at infinity. Distributed along the unit
circle is the 10-point orbit of 72-points. The 5 lines of reflective
symmetry passing through (0,0) correspond to the 5 \RP{1}\ intersections
of \RR\ and the 45-lines whose associated involution preserves \RR.

\newpage


\begin{figure}[hp]

\epsfxsize=\figWidth
 \epsfbox{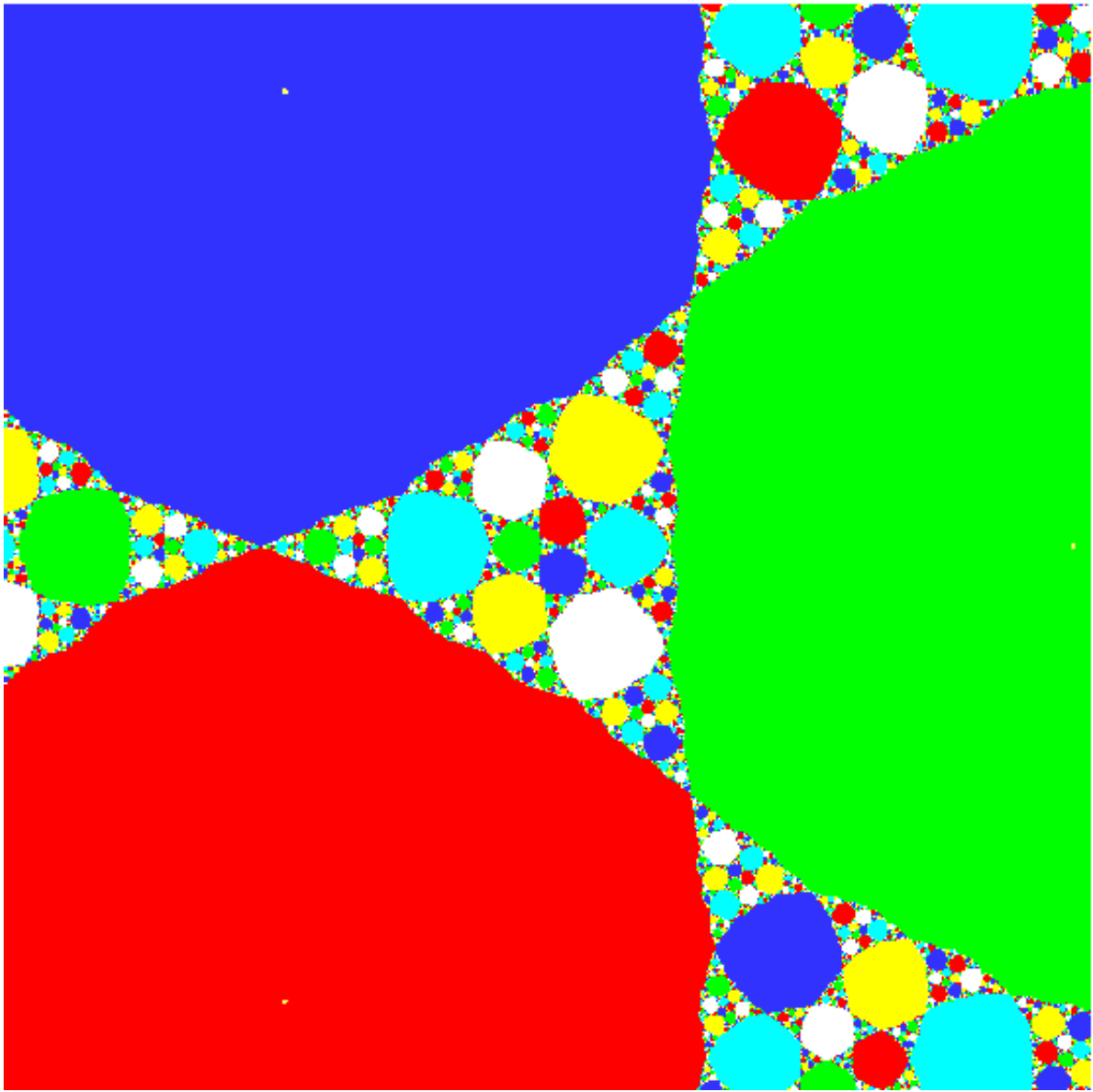}

\caption{Icosahedral dynamics}

\label{fig:h19Conic}

When restricted to an icosahedral conic, the canonical degree 19 map
$h_{19}$ has icosahedral symmetry and attracts almost all points on the
conic to antipodal pairs of vertices. Each of the 6 colors corresponds
to such a pair and the 3 large basins are immediate---each contains a
superattracting vertex.
 
\end{figure}

\newpage


\begin{figure}[hp]

\epsfxsize=\figWidth
 \epsfbox{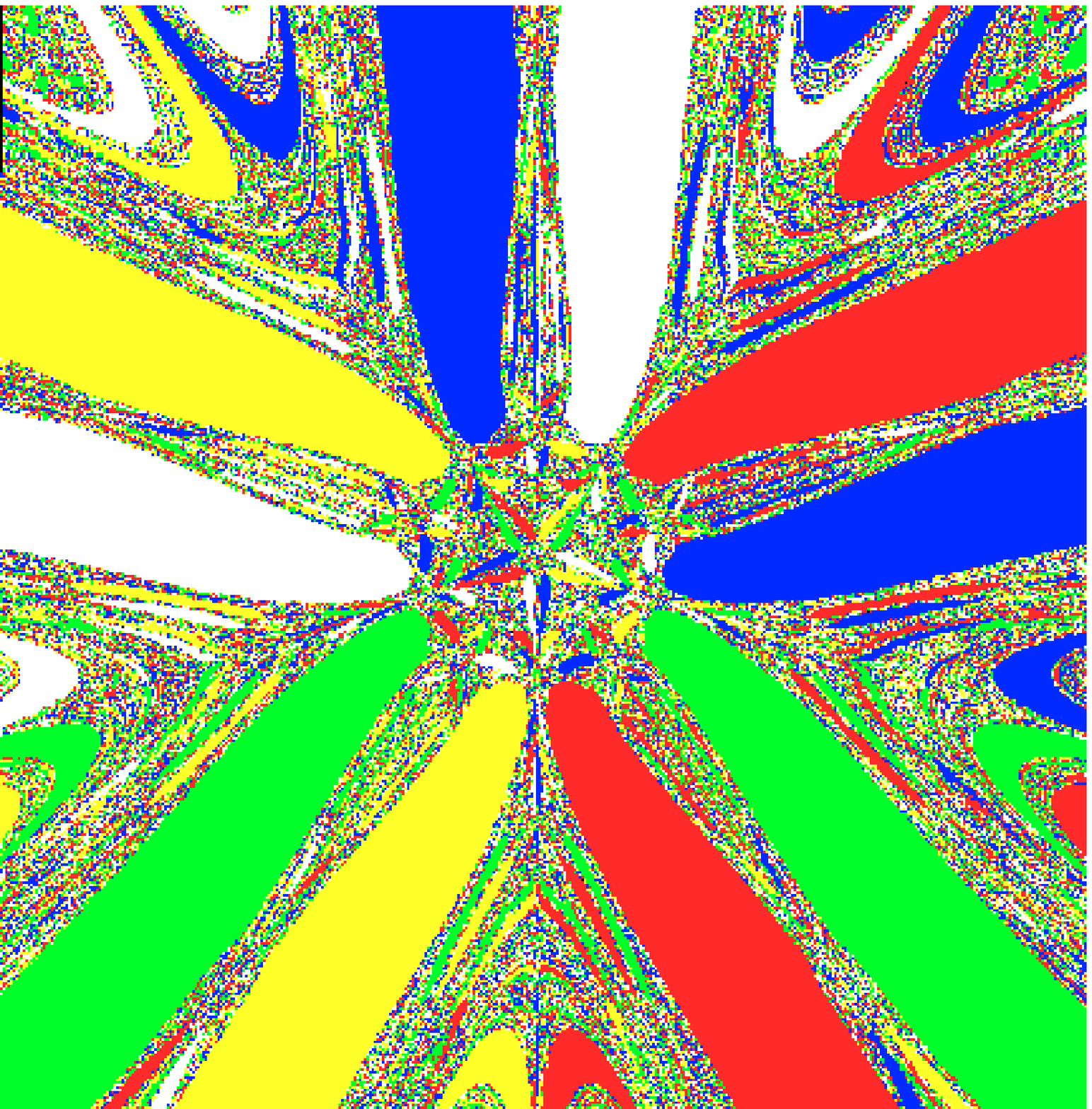}

\caption{\RP{2}\ dynamics}

\label{fig:h19RP2}

For this plot the vertical and horizontal scales are roughly
from -2 to 2.  The large ``radial" basins are immediate, i.e., each
contains 1 of the 72-points and come in pairs as do the period 2
attractors.  Notice the repulsive behavior along the 45-lines and
particularly at their 5-fold intersection at the 36-point.

\end{figure}